\documentclass[11pt,tgrind,psbox]{article}
\usepackage{latexsym}
\usepackage{graphics}
\usepackage{amsmath}
\usepackage{xspace}
\usepackage{amssymb}
\usepackage{psfrag}
\usepackage{epsfig}
\usepackage{fullpage}
\usepackage{amsmath,amsthm}
\usepackage{epsfig}
\usepackage{pst-all}

\newcommand {\bel}[1]{\begin{align*}}
\newcommand {\eel}[1]{\end{align*}}
\newcommand {\bea}{\begin{eqnarray}}
\newcommand {\eea}{\end{eqnarray}}

\newcommand{\pr}{\mathbb{P}}

\newcommand{\G}{\mathbb{G}}
\newcommand{\T}{\mathbb{T}}
\newcommand{\R}{\mathbb{R}}

\newcommand{\M}{{\bf {\mathbb{M}}}}
\newcommand{\mb}[1]{\mbox{\boldmath $#1$}}

\newcommand{\mbw}{\mb{w}}
\newcommand{\mbx}{\mb{x}}
\newcommand{\mbX}{\mb{X}}
\newcommand{\mby}{\mb{y}}
\newcommand{\mbz}{\mb{z}}
\newcommand{\mbL}{\mb{L}}

\newcommand{\X}{\mathcal{X}}

\newcommand{\ignore}[1]{\relax}

\newtheorem{theorem}{Theorem}

\newtheorem{lemma}{Lemma}

\newtheorem{prop}{Proposition}
\newtheorem{coro}{Corollary}

\newtheorem{Defi}{Definition}

\title{Correlation decay and deterministic FPTAS for counting list-colorings of a graph
\footnote{A Preliminary version of this paper appeared in
Proceedings of 18th ACM-SIAM Symposium on Discrete Algorithms (SODA)~\cite{GamarnikKatz}}}

\author{
 {\sf David Gamarnik }
  \thanks{Operations Research Center and Sloan School of Management, MIT, Cambridge, MA,  02139, e-mail: {\tt
gamarnik@mit.edu}}\and
{\sf Dmitriy Katz} \thanks{Operations Research Center, MIT, Cambridge, MA,  02139, e-mail: {\tt
dimdim@mit.edu}}}

\begin{document}

\maketitle

\begin{abstract}
We propose a deterministic algorithm for approximately counting the
number of list colorings of a graph. Under the assumption that the
graph is triangle free, the size of every list is at least $\alpha
\Delta$, where $\alpha$ is an arbitrary constant bigger than
$\alpha^{**}=2.8432\ldots$, and $\Delta$ is the maximum degree of the graph, we
obtain the following results. For the case when the size of the each
list is a large  constant, we show the existence of a
\emph{deterministic} FPTAS for computing the total number of list
colorings. The same deterministic algorithm has complexity
$2^{O(\log^2 n)}$, without any assumptions on the sizes of the
lists, where $n$ is the instance size. We further extend our method
to a discrete Markov random field (MRF) model. Under certain assumptions
relating the size of the alphabet, the degree of the graph and the interacting potential
we again construct a deterministic FPTAS for computing the partition function
of a  MRF.

Our results are not based on the most powerful existing counting
technique -- rapidly mixing Markov chain method. Rather we build
upon concepts from statistical physics, in particular, the decay of
correlation phenomena and its implication for the uniqueness of
Gibbs measures in infinite graphs. This approach was proposed in two
recent papers \cite{BandyopadhyayGamarnikCounting} and
\cite{weitzCounting}. The principle insight of this approach is
that the correlation decay property can be established with respect
to certain \emph{computation tree}, as opposed to the conventional
correlation decay property with
respect to graph theoretic neighborhoods of a  given node. This
allows truncation of  computation at a logarithmic depth in order to
obtain polynomial accuracy in polynomial time.
\end{abstract}


\section{Introduction}
This paper is devoted to the problem of computing the total number of list colorings of a graph.
It is further concerned with the problem of computing a partition function
corresponding to a Markov random field  (also known as graphical) model.
The setting for the list coloring problem is as follows. Each node of a given  graph is associated with a list of colors.
An assignment of nodes to colors is called list coloring if every node is assigned to some color
from its list and no two nodes sharing an edge are assigned to the same color.
When all the lists are identical, the problem reduces to the problem of coloring of a graph.
The problem of determining
whether a list coloring exists is NP-hard, but provided that the size of each list is stictly larger
than the degree for each node, a simple greedy algorithm produces a coloring.
We are concerned with the corresponding counting problem -- compute
the total number of list colorings of a given graph/list pair. This problem is known to be $\#P$ hard
even for the restricted problem of counting the colorings,
and the focus is on the approximation algorithms. The existing approximation schemes
are mostly based on the rapidly mixing Markov chain
technique, also known as Glauber dynamics approach. It was established by Jerrum~\cite{Jerrum} that
the Glauber dynamics corresponding to graphs where the ratio of the number of colors to degree satisfies
$q/\Delta\geq 2$, mixes rapidly. This leads to a randomized approximation algorithm for enumerating the number of colorings.
The $2$-barrier was first broken by Vigoda~\cite{Vigoda}, who lowered the ratio requirement to $11/6$. Many further significant
improvements were obtained subsequently. The state of the art is summarized in \cite{FriezeVigodaSurvey}.
For a while the improvement over $11/6$ ratio came at a cost of lower bound $\Omega(\log n)$ on the maximum degree,
where $n$ is the number of nodes. This requirement was lifted by Dyer et al. \cite{DyerFriezeHayesVidoga}.

In this paper we focus on a different approach to the  counting list colorings problem.
Our setting is a list coloring problem. We require that the size of every list is at least
$\alpha\Delta+\beta$, where $\alpha>\alpha^{**}=2.8432\ldots$ - the unique solution to $\alpha e^{-{1\over \alpha}}=2$,
and $\beta$ is a large constant which depends on $\alpha-\alpha^{**}$.
Our girth restriction is   $g\geq 4$, namely, the graph is triangle-free.
We obtain the following results. First, assuming that the size of each list is at most a constant, we construct a
\emph{deterministic} Fully Polynomial Time Approximation Scheme (FPTAS) for the problem of computing the total
number of list colorings of a given graph/list pair. Second, for arbitrary graph/list pair (no
assumptions on the list sizes) we construct an approximation algorithm with complexity
$2^{O(\log^2 n)}$. Namely, our algorithm is super-polynomial but still significantly quicker than exponential time.

Although our regime $\alpha>2.8432\ldots$ is weaker than $q/\Delta>2$, for which the Markov chain
is  known  to mix rapidly, the important contribution of our method is that it provides a deterministic algorithm.
Presently no deterministic algorithms are known for counting approximately the number of coloring of a graph.

Our approach is based on establishing a certain \emph{correlation decay} property which has been considered in many settings
\cite{BrightwellWinklerHomomorphisms2004},
including the coloring problem \cite{SalasSokal}, \cite{GoldbergMartinPaterson},
\cite{BrightwellWinklerColoring},\cite{JonassonColoring2002}
and has been recently a subject of interest. In particular, the correlation decay has
been established in \cite{GoldbergMartinPaterson} for coloring  triangle-free graphs under
the assumption that $\alpha>\alpha^*=1.763\ldots$, the unique solution of $\alpha e^{-{1\over \alpha}}=1$.
(Some mild additional assumptions were adopted).
The principal motivation for establishing the
correlation decay property comes from statistical physics, in particular
the connection with the uniqueness of the associated Gibbs measure
(uniform measure in our setting) on infinite
versions of the graph, typically lattices. Recently, however,  a new approach linking
correlation decay to counting algorithms was proposed
in Bandyopadhyay and Gamarnik~\cite{BandyopadhyayGamarnikCounting} and Weitz~\cite{weitzCounting}. The idea is
to use correlation decay property instead of Markov sampling for computing marginals of the Gibbs (uniform) distribution.
This leads to a deterministic approach since the marginals are computed using a dynamic programming like scheme
(also known as Belief Propagation (BP) algorithm \cite{YedidiaFreemanWeiss}). This approach typically needs a locally-tree
like structure (large girth) \cite{ShahRandomMatching} in order to
be successful. The large girth assumption was explicitly assumed in \cite{BandyopadhyayGamarnikCounting}, where the problems
of computing the number of independent sets and colorings in some special structured (regular) graphs was considered.
Weitz~\cite{weitzCounting} cleverly by-passes the large girth assumption by using a certain self-avoiding tree construction
thus essentially reducing the problem to a problem on a tree with careful boundary conditions implied by independent sets.
This idea was used recently by Jung and Shah~\cite{JungShah} to introduce a version of a
BP algorithm which works on a non-locally-tree like
graphs, where appropriate correlation decay can be established.
This approach works for binary type problems (independent sets, matchings, Ising model) but does not apparently extend
to multi-valued problems.

In this paper we propose a general  deterministic approximate counting algorithm which can be used for arbitrary multi-valued
counting problem. We also by-pass the large
girth assumption by considering a certain \emph{computation tree} corresponding to the Gibbs (uniform for the case of colorings)
measure. Our principal insight is establishing correlation decay for the computation tree
as opposed to the conventional correlation decay associated with the
graph-theoretic structure of the graph. We provide a discussion explaining why it is crucial to
establish the correlation decay in this
way in order to obtain FPTAS. Contrast
this with \cite{GoldbergMartinPaterson} where correlation decay is established for the coloring
problem but in the conventional graph-theoretic distance
sense. Our method is similar to the self-avoiding walk method of Weitz's but somewhat more direct as
the step of relating the marginal probability on a graph to the marginal probability on the tree is by-passed
in our computation tree approach.
The advantage of establishing correlation decay on a computation tree as opposed to the original graph has been highlighted also
in \cite{JordanTatikonda} in the context of BP algorithms and the Dobrushin's Uniqueness condition.
More importantly our approach works for general, not necessarily two-valued model. We extend our approach
to Markov random field model and also show that under some conditions the computation tree satisfies the
correlation decay property and, as a result, one obtains a deterministic algorithm for computing
approximately the associated partition function.

The remainder of the paper has the following structure. The model description and the main result are stated
in Section~\ref{section:MainResult}. Some preliminary technical results
are established in Section~\ref{section:Preliminary}. The description of the algorithm and its complexity are
subject of Section~\ref{section:AlgorithmComplexity}. The principal technical result is established in
Section~\ref{section:CD}. The key result is Theorem~\ref{theorem:correlationdecay}, which establishes
the correlation decay result on a computation tree arising in computing the marginals of the uniform distribution
on the set of all list colorings.
Section~\ref{section:MRF} is devoted to the extension of our approach to a Markov random field.
Section~\ref{section:comparison} provides a brief comparison between the correlation decay on a computation tree
and the correlation decay in a conventional sense.
Some conclusions and open problems are in Section~\ref{section:conclusions}.

\section{Definitions and the main result}\label{section:MainResult}
We consider a simple graph $\mathbb{G}$ with the node set $V=\{v_1,v_2,\ldots,v_{|V|}\}$.
Our graph is assumed to be triangle-free. Namely the girth (the size of the smallest cycle) is at least $g\geq 4$.
Let $E,\Delta$ denote respectively
the set of  edges and the maximum degree of the graph.  $\Delta(v)$ denotes the degree
of the node $v$.
Each node $v$ is associated with a list of colors
$L(v)\subset \{1,2,\ldots,q\}=\cup_{v\in V}L(v)$, where $\{1,2,\ldots,q\}$ is the total universe of colors.
We let $\mb{L}=(L(v), 1\leq v\leq n)$ denote the vector of lists.
We also let $\|\mbL\|=\max_v |L(v)|$ the size of the largest list.
The list-coloring problem on $\G$ is formulated
as follows:  associate each node $v$ with a color
$c(v)\in L(v)$ such that no two nodes sharing an edge are associated
with the same color. When all the lists are identical and contain $q$ elements,
the corresponding problem is the problem of coloring $\G$ using $q$ colors.
We let $|L(v)|$ denote the cardinality of $L(v)$.
It is easy to see that if
\begin{align}\label{eq:degreePlusone}
|L(v)|\geq \Delta(v)+1
\end{align}
for every node $v$,
then a simple greedy procedure
produces a list-coloring. We adopt here a stronger assumption
\begin{align}\label{eq:degreealpha}
|L(v)|\geq \alpha \Delta(v)+\beta,
\end{align}
where $\alpha$ is an arbitrary constant strictly larger than $\alpha^{**}$, the unique solution of
$\alpha^{**}\exp(-{1\over \alpha*})=2$. That is $\alpha^{**}\approx 2.8432\ldots$. We also assume that
$\beta$ is a large constant which depends on $\alpha$.
To be more specific we assume that $\beta=\beta(\alpha)$ is large enough to satisfy
\begin{align}\label{eq:beta}
(1-{1\over \beta})\alpha e^{-{1\over \alpha}(1+{1\over \beta})}>2,
\end{align}
which is always possible when $\alpha>\alpha^{**}$.

Let $Z(\G,\mb{L})$ denote the total number of possible list-colorings of a graph/list pair
$(\G,\mbL)$. The corresponding counting problem is to compute
(approximately) $Z(\G,\mb{L})$. In statistical physics terminology, $Z(\G,\mb{L})$ is the partition function.
 We let $Z(\G,\mb{L},\chi)$ denote the number of list
colorings of $(\G,\mbL)$ which satisfy some condition $\chi$. For example $Z(\G,\mb{L},c(v)=i,c(u)=j)$
is the number of list colorings such that the color of $v$ is $i$  and the color of $u$ is $j$.

On the space of all list colorings of $\G$ we consider a uniform probability distribution, where each
list coloring assumes weight $1/Z(\G,\mb{L})$. For every node/color pair $v\in V, i\in L(v)$,
$\pr_{\G,\mbL}(c(v)=i)$ denotes the probability that node $v$ is colored $i$ with respect
to this probability measure.
The size of the instance corresponding to a graph/list pair$(\G,\mbL)$ is defined to be
$n=\max\{|V|,|E|,q\}$.

\begin{Defi}
An approximation algorithm ${\cal A}$ is defined to be a Fully Polynomial Time Approximation Scheme for a computing
$Z(\G,\mbL)$ if given arbitrary $\delta>0$ it produces a value $\hat Z$ satisfying
\begin{align*}
1-\delta\leq {\hat Z\over Z(\G,\mbL)}\leq 1+\delta,
\end{align*}
in time which is polynomial in $n,{1\over \delta}$.
\end{Defi}

We now state our main result.
\begin{theorem}\label{theorem:MainResult}
There exist a deterministic algorithm which provides a FPTAS for computing $Z(\G,\mb{L})$ for arbitrary
graph list pair $\G,\mbL$ satisfying (\ref{eq:degreealpha}), when the size of the largest list  $\|\mbL\|$ is
constant. The same algorithm has complexity $2^{O(\log^2n)}$, without any  restriction  on $\|\mbL\|$,
 where $n$ is the size of the instance.
\end{theorem}

\section{Preliminary technical results}\label{section:Preliminary}

\subsection{Basic recursion}\label{subsection:recursion}
We begin by establishing a standard relationship between the partition function $Z(\G,\mbL)$ and the marginals
$\pr_{\G,\mbL}(c(v)=i)$. The relation, also known as cavity method, is also the basis of the Glauber dynamics
approach for computing partition functions.

\begin{prop}\label{prop:cavity}
Consider an arbitrary list coloring  $i_1,\ldots,i_{|V|}$ of the graph $\G$ (which can be constructed
using a simple greedy procedure).
For every $k=0,1,\ldots,|V|-1$ consider a graph list pair $\G_k,\mbL_k$, where
$(\G_0,\mbL_0)=(\G,\mbL)$,
$\G_k=\G\setminus\{v_1,\ldots,v_k\}, k\geq 1$ and
the list $\mbL_k$ is obtained by deleting from each list $L(v_l),l>k$ a color $i_r,r\leq k$ if $(v_l,v_r)\in E$.
Then
\begin{align*}
Z(\G,\mbL)=\prod_{0\leq k\leq |V|-1}\pr_{\G_k,\mbL_k}^{-1}(c(v_k)=i_k).
\end{align*}
\end{prop}

\begin{proof}
We have
\begin{align*}
\pr_{\G,\mbL}(c(v_1)=i_1)={Z(\G,\mbL,c(v_1)=i_1)\over Z(\G,\mbL)}={Z(\G_1,\mbL_1)\over Z(\G,\mbL)},
\end{align*}
from which we obtain
\begin{align*}
Z(\G,\mbL)=\pr_{\G,\mbL}(c(v_1)=i_1)^{-1}Z(\G_1,\mbL_1).
\end{align*}
Iterating further for $k\geq 2$ we obtain the result.
\end{proof}

Our algorithm is based on a recursive procedure which relates the number of list colorings of a given graph/list pair
in terms of the number of list colorings of some reduced graph/list pairs.

Given a pair $(\G,\mbL)$ and a node $v\in \G$, let $v_1,\ldots,v_m$ be the set of neighbors of $v$.
For every pair $(k,i)\in \{1,\ldots,m\}\times L(v)$ we define a new pair $(\G_v,\mbL_{k,i})$ as follows.
The set of nodes of $\G$ is
$V_{k}=V\setminus \{v\}$ and  $L_{k,i}(v_r)=L(v_r)\setminus \{i\}$ for $1\leq r<k$,
 $L_{k,j}(u)=L(u)$ for all other $u$.
Namely, we first delete node $v$ from the graph. Then we delete color $i$ from the lists corresponding to the nodes
$v_r,r<k$, and leave all the other lists intact.
\begin{lemma}\label{lemma:degreealphareduced}
The graph/list pair $(\G_v,\mbL_{k,j})$  satisfies (\ref{eq:degreealpha}) for every $1\leq k\leq m,j\in L(v)$,
provided that $(\G,\mbL)$ does.
\end{lemma}

\begin{proof}
When we create  graph $\G_v$ from $\G$ the list size of every remaining node
either stays the same or is reduced by one. The second event can only happen for
neighbors $v_1,\ldots,v_m$ of the deleted node $v$.
When the list is reduced by one the degree is reduced by one as well. Since $\alpha>1$, the assertion follows
by observing that $|L(v_k)|\geq \alpha \Delta(v_k)+\beta$ implies $|L(v_k)|-1\geq \alpha (\Delta(v_k)-1)+\beta$.
\end{proof}

The  basis of our algorithm is the following simple result.
\begin{prop}\label{prop:Recursion1}
Given a  graph/list pair $(\G,\mbL)$ and a  node $v$, suppose $\Delta(v)=m>0$.
For every  $i\in L(v)$
\begin{align}\label{eq:Aalpha}
\pr_{\G,\mbL}(c(v)=i)={\prod_{1\leq k\leq m}(1-\pr_{\G_v,\mbL_{k,i}}(c(v)=i))\over
\sum_{j\in L(v)}\prod_{1\leq k\leq m}(1-\pr_{\G_v,\mbL_{k,j}}(c(v)=j))}.
\end{align}
\end{prop}
The recursion as well as the proof is similar to the one used by Weitz in \cite{weitzCounting}, except we bypass the construction
of a self-avoiding tree, considered in \cite{weitzCounting}.

\begin{proof}
Consider a graph/list $(\G_v,\mbL)$ obtained simply by removing node $v$ from $\G$, and leaving $\mbL$ intact for the
remaining nodes.
We have
\begin{align*}
\pr_{\G,\mbL}(c(v)=i)&={\pr_{\G,\mbL}(c(v)=i)\over \sum_{j\in L(v)}\pr_{\G,\mbL}(c(v)=j)} \\
&={Z(\G,\mbL,c(v)=i)Z^{-1}(\G,\mbL)\over \sum_{j\in L(v)}Z(\G,\mbL,c(v)=j)Z^{-1}(\G,\mbL)} \\
&={Z(\G_v,\mbL,c(v_k)\neq i, ~1\leq k\leq m)\over \sum_{j\in L(v)}Z(\G_v,\mbL,c(v_k)\neq j, ~1\leq k\leq m)} \\
&={\pr_{\G_v,\mbL}(c(v_k)\neq i, ~1\leq k\leq m)\over \sum_{j\in L(v)}\pr_{\G_v,\mbL}(c(v_k)\neq j, ~1\leq k\leq m)}
\end{align*}
Now, for every $j\in L(v)$
\begin{align*}
\pr_{\G_v,\mbL}(c(v_k)\neq j, ~1\leq k\leq m)=
\pr_{\G_v,\mbL}(c(v_1)\neq j)\prod_{2\leq k\leq m}\pr_{\G_v,\mbL}(c(v_k)\neq j|c(v_r)\neq j, ~1\leq r<k)
\end{align*}
We observe that $\mbL_{1,j}=\mbL$ for every $j$ (no colors are removed due
to the vacuous condition $r<1$), and
$\pr_{\G_v,\mbL}(c(v_k)\neq j|c(v_r)\neq j, ~1\leq r<k)=\pr_{\G_v,\mbL_{k,j}}(c(v_k)\neq j)$.
Namely
\begin{align*}
\pr_{\G_v,\mbL}(c(v_k)\neq j, ~1\leq k\leq m)=
\prod_{1\leq k\leq m}\pr_{\G_v,\mbL_{k,j}}(c(v_k)\neq j)=
\prod_{1\leq k\leq m}(1-\pr_{\G_v,\mbL_{k,j}}(c(v_k)=j)).
\end{align*}
Substituting this expression we complete the proof.
\end{proof}

\subsection{Upper and lower bounds}
The condition (\ref{eq:degreealpha}) allows us to obtain the following simple bounds.
\begin{lemma}\label{lemma:lessbeta}
For every $\G,\mbL$, node $v$ and a color $i\in L(v)$
\begin{align*}
\pr_{\G,\mbL}(c(v)=i)\leq {1\over \beta}.
\end{align*}
\end{lemma}

\begin{proof}
Observe that given an arbitrary coloring of the neighbors $v_1,\ldots,v_m$ of $v$,
there are at least $|L(v)|-\Delta(v)\geq \beta$ colors remaining. Then the upper bound holds.
\end{proof}

From this simple bound we now establish a different upper bound and also a lower bound using the triangle free assumption.
\begin{lemma}\label{lemma:lessdelta}
There exist $\epsilon_0=\epsilon_0(\alpha)\in (0,1)$ and $\beta>0$ such that
for every $\G,\mbL$, node $v$ and a color $i\in L(v)$
\begin{align*}
q^{-1}(1-\beta^{-1})^{\Delta}\leq \pr_{\G,\mbL}(c(v)=i)\leq {1\over 2\Delta(v)(1+\epsilon_0)}.
\end{align*}
\end{lemma}
We note that the upper bounds of this lemma and Lemma~\ref{lemma:lessbeta} are not comparable, since values
of $\Delta(v)$ could be smaller and larger than $\beta$.

\begin{proof}
We let $v_1,\ldots,v_m$ denote the neighbors of $v$, $m=\Delta(v)$ and let
$v_{kr}$ denote the set of neighbors of $v_k$, other than $v$ for $k=1,\ldots,m$. We will establish that
for any coloring of nodes $(v_{kr})$, which we generically denote by $\mb{c}$, we have
\begin{align*}
q^{-1}(1-\beta^{-1})^{\Delta}\leq \pr_{\G,\mbL}(c(v)=i|\mb{c})\leq {1\over 2m(1+\epsilon_0)}.
\end{align*}
The corresponding inequality for the unconditional
probability then follows immediately. Now observe that, since the girth is at least $4$, then there are no edges between
$v_k$. Then $\pr_{\G,\mbL}(c(v)=i|\mb{c})$ is the probability $\pr_{\T}(c(v)=i)$ that $v$ is colored
$i$ in a depth-1 tree $\T\triangleq \{v,v_1,\ldots,v_m\}$, where the lists $\hat L(v_k)$
of $v_k$ are obtained  from $L(v_k)$ by deleting
the colors used by the neighbors $v_{kr}$ by coloring $\mb{c}$. From the assumption (\ref{eq:degreealpha})
we have that the remaining lists $\hat L(v_k)$ have size at least $|L(v_k)|-\Delta(v_k)\geq \beta$ each.
Let $t_i=\pr_{\T}(c(v)=i)$.
For each color $j\in L(v)$ let
$t_{j,k}=1/|\hat L(v_k)|$ if $j\in \hat L(v_k)$ and $=0$ otherwise. Proposition~\ref{prop:Recursion1}
then simplifies to
\begin{align}\label{eq:boundti}
t_i={\prod_{1\leq k\leq m} (1-t_{i,k})\over \sum_{j\in L(v)}\prod_{1\leq k\leq m}
(1-t_{j,k})}\leq {1\over \sum_{j\in L(v)}\prod_{1\leq k\leq m} (1-t_{j,k})},
\end{align}
for every $i\in L(v)$, where $\prod_{1\leq k\leq m}$ is defined to be equal to unity when $m=0$.
From the equality part, applying $t_{j,k}\leq 1/\beta$, we get
\begin{align*}
t_i\geq |L(v)|^{-1}(1-\beta^{-1})^m\geq q^{-1}(1-\beta^{-1})^{\Delta},
\end{align*}
and the lower bound is established.

We now focus on the upper bound and use the inequality part of (\ref{eq:boundti}).
Thus it suffices to show that
\begin{align}
\sum_{j\in L(v)}\prod_k (1-t_{j,k})\geq 2(1+\epsilon_0)m \label{eq:geq2epsilon0}
\end{align}
for some constant $\epsilon_0>0$.
Using the first order Taylor expansion for $\log z$ around $z=1$,
\begin{align*}
\prod_{1\leq k\leq m} (1-t_{j,k})&=\prod_{1\leq k\leq m}e^{\log (1-t_{j,k})} \\
&= \prod_{1\leq k\leq m}e^{-t_{j,k}-{1\over 2(1-\theta_{j,k})^2}t_{j,k}^2},
\end{align*}
for some $0\leq \theta_{j,k}\leq t_{j,k}$, since $-1/z^2$ is the second derivative of $\log z$.
Again using the bound $t_{j,k}\leq 1/\beta$, we have
$(1-\theta_{j,k})^2 \geq (1-1/\beta)^2$.    
We assume that $\beta$ is a sufficiently large constant ensuring $(1-1/\beta)^2>1/2$.
Thus we obtain the following lower bound
\begin{align*}
\prod_{1\leq k\leq m} (1-t_{j,k})\geq \prod_{1\leq k\leq m}e^{-t_{j,k}-{t_{j,k}^2\over 2(1-1/\beta)^2}}
\geq e^{-(1+{1\over \beta})\sum_k t_{j,k}}
\triangleq e^{-(1+{1\over \beta})T_{j}},
\end{align*}
where $T_{j}$ stands for $\sum_k t_{j,k}$. Then
\begin{align*}
\sum_{j\in L(v)}\prod_{1\leq k\leq m} (1-t_{j,k})\geq
\sum_{j\in L(v)}e^{-(1+{1\over \beta})T_{j}}\geq |L(v)| e^{-{1\over |L(v)|}(1+{1\over \beta})\sum_jT_{j}},
\end{align*}
where we have used an inequality between the average arithmetic and average geometric.
Finally  we observe
\begin{align*}
\sum_{j\in L(v)} T_j=\sum_{j,k}t_{j,k}=\sum_{1\leq k\leq m}\sum_{j\in \hat L(v_k)}{1\over |\hat L(v_k)|}= m.
\end{align*}

Thus
\begin{align*}
\sum_{j\in L(v)}\prod_{1\leq k\leq m} (1-t_{j,k})
\geq |L(v)| e^{-{m\over |L(v)|}(1+{1\over \beta})}\geq (\alpha m+\beta)e^{-{1\over \alpha}(1+{1\over \beta})}
>\alpha me^{-{1\over \alpha}(1+{1\over \beta})}
\end{align*}
The condition $\alpha>\alpha^{**}$ implies that there exists  a sufficiently large $\beta$
such that $\alpha e^{-{1\over \alpha}(1+{1\over \beta})}>2$. We find $0<\epsilon_0<.1$
such that $\alpha e^{-{1\over \alpha}(1+{1\over \beta})}=2(1+\epsilon_0)$.
We obtain a required lower bound (\ref{eq:geq2epsilon0}).
\end{proof}

\section{Algorithm and complexity}\label{section:AlgorithmComplexity}

\subsection{Description of an algorithm}

Our algorithm is based on the idea of trying to approximate the value of $\pr_{\G,\mbL}(c(v)=i)$,
by performing a certain recursive
computation using (\ref{eq:Aalpha}) a fixed number of times $d$ and then using a correlation decay principle to guarantee
the accuracy of the approximation. Specifically, introduce
a function $\Phi$ which takes as an input a vector $(\G,\mbL,v,i,d)$
and takes some values
$\Phi(\G,\mbL,v,i,d)\in [0,1]$. The input $(\G,\mbL,v,i,d)$ to $\Phi$ is any vector, such that
such that $v$ is a node in $\G$,
$i$ is an arbitrary color, and $d$ is an arbitrary non-negative integer. Function $\Phi$ is defined recursively
in $d$. The quantity $\Phi$ "attempts" to approximate
$\pr_{\G,\mbL}(c(v)=i)$. The quality of the approximation is controlled by $d$.
We define $\Phi$ as follows. For every input $(\G,\mbL,v,i,d)$ such that $i\notin L(v)$ we set
$\Phi(\G,\mbL,v,i,d)=0$. Otherwise we set the values as follows.

\begin{itemize}
\item When $d=0$, we set
$\Phi(\G,\mbL,v,i,d)=1/|L(v)|$ for every input $(\G,\mbL,v,i)$. (It turns out that for our application
the initialization values are not important, due to the decay of correlations).

\item For every $d\geq 1$, if $\Delta(v)=0$, then $\Phi(\G,\mbL,v,i,d)=1/|L(v)|$ for all $i\in L(v)$.
Suppose $\Delta(v)=m>0$
and $v_1,\ldots,v_m$ are the neighbors of $v$. Then for every $i\in L(v)$ we define
\begin{align}\label{eq:definePhi}
\Phi(\G,\mbL,v,i,d)=\min\Big[{1\over 2(1+\epsilon_0)m},{1\over \beta},{\prod_{1\leq k\leq m}
(1-\Phi(\G_v,\mbL_{k,i},v_k,i,d-1))\over
\sum_{j\in L(v)}\prod_{1\leq k\leq m}(1-\Phi(\G_v,\mbL_{k,j},v_k,j,d-1))}\Big].
\end{align}
The last part of the expression inside $\min[\cdot]$ corresponds directly to the expression
(\ref{eq:Aalpha}) of Proposition~\ref{prop:Recursion1}. Specifically, if it was true
that $\Phi(\G_v,\mbL_{k,j},v_k,j,d-1)=\pr_{\G_v,\mbL_{k,j}}(c(v_k)=j)$,
then, by Lemmas~\ref{lemma:lessbeta},\ref{lemma:lessdelta},
the minimum in (\ref{eq:definePhi}) would be achieved by the third expression,
and  then the value of
$\Phi(\G,\mbL,v,i,d)$ would be exactly $\pr_{\G,\mbL}(c(v)=i)$.
\end{itemize}

We will use the correlation decay
property to establish that the difference between the two values, modulo rescaling, is diminishing
as $d\rightarrow\infty$. Note that the computation of $\Phi$ can be done recursively in $d$ and it involves
a dynamic programming type recursion. The underlying computation is done essentially on a tree of graph
list pairs $\G_s,\mbL_s$ generated during the recursion. We refer to this tree as
\emph{computation tree} with depth $d$.


We now describe our algorithm for approximately computing $Z(\G,\mbL)$.
The algorithm is parametrized by the "quality" parameter $d$.

\vspace{.1in}

\textbf{Algorithm CountCOLOR}
\vspace{.1in}

{\tt
INPUT: A graph/list pair $(\G,\mbL)$ and a positive integer $d$.

BEGIN

Set $\hat Z=1, \hat \G=\G,\hat \mbL=\mbL$.

While $\hat G\neq \emptyset$, find an arbitrary node $v\in \hat G$ and a color $i\in \hat L(v)$. Compute
\begin{align}\label{eq:pvfromPhi}
\hat p(v,i)\triangleq \Phi(\hat G,\hat\mbL,v,i,d).    
\end{align}
Set $\hat Z=\hat p^{-1}(v,i)\hat Z, \hat G=\hat G\setminus \{v\},
\hat L(u)=\hat L(u)\setminus \{i\}$ for all neighbors $u$ of $v$ in $\hat G$, and $\hat L(u)$ remains the same for all
other nodes.

END

OUTPUT: $\hat Z$.}

\vspace{.2in}

\subsection{Some properties}
We now establish some properties of $\Phi$.
\begin{lemma}\label{lemma:Phiproperties}
The following holds for every $\G,\mbL,v,i\in L(v), d\geq 0$.
\begin{align}
\Phi(\G,\mbL,v,i,d)&\leq \min\big[{1\over \beta},{1\over 2(1+\epsilon_0)\Delta(v)}\big], \label{eq:PhiProperties1}\\
\sum_{i\in L(v)}\Phi(\G,\mbL,v,i,d)&\leq 1, \label{eq:PhiProperties2}\\
\Phi(\G,\mbL,v,i,d)&\geq q^{-1}(1-1/\beta)^\Delta. \label{eq:lowerboundPhi}
\end{align}
\end{lemma}

\begin{proof}
(\ref{eq:PhiProperties2}) follows directly from the definition of $\Phi$.
To show (\ref{eq:PhiProperties1})
we consider cases.
For $d\geq 1$ this follow directly from the recursion (\ref{eq:definePhi}). For $d=0$,
this follows since $\Phi(\G,\mbL,v,i,0)=1/|L(v)|\leq 1/(\alpha\Delta(v)+\beta)$
and $2(1+\epsilon_0)<2.2<\alpha$.
We now establish (\ref{eq:lowerboundPhi}).
For the case $d=0$ this follows since $1/|L(v)|\geq 1/q$. For the case $d\geq 1$ this follows from the
recursion (\ref{eq:definePhi}) since $1/\beta,1/(2(1+\epsilon_0)\Delta(v))>1/q$ and the third term
inside the minimum operator is at least $q^{-1}(1-1/\beta)^\Delta$, using upper bound
$\Phi(\G,\mbL,v,i,d-1)\leq 1/\beta$ which we have from (\ref{eq:PhiProperties1}).

\end{proof}

\subsection{Complexity}
We begin by analyzing the complexity of computing function $\Phi$. Recall that  $n=\max(|V|,|E|,q)$ is the size
of the instance.

\begin{prop}
For any given node $v$, the function $\Phi$ can be computed in time $2^{O(d(\log \|L\|+\log\Delta))}$.
In particular when $d=O(\log n)$, the overall computation is $2^{O(\log^2n)}$.
If in addition the size of the largest list  $\|L\|$ is constant then the computation time
is polynomial in $n$.
\end{prop}

\begin{proof}
Let $T(d)$ denote the complexity of computing function $\Phi(\cdot,d)$. Clearly,
$T(0)=O(\|L\|)$. We now express $T(d)$ in terms of $T(d-1)$. Given a node $v$,
in order to compute $\Phi(\G,\mbL,v,i,d)$ we first identify the neighbors $v_1,\ldots,v_m$ of $v$.
Then we create graph/list pairs $\G_v,\mbL_{j,k}, 1\leq k\leq m, j\in L(v)$,   compute
$\Phi(\cdot,d-1)$ for each of this graphs, and use this to compute $\Phi(\G,\mbL,v,i,d)$.
The overall computation effort is then
\begin{align*}
T(d)&=O(\|L\|\Delta T(d-1)).  
\end{align*}
Iterating over $d$ we obtain $T(d)=O(\|L\|^{d+1}\Delta^d)=O(2^{(d+1)(\log \|L\|+\log\Delta)})=2^{O(d(\log \|L\|+\Delta))}$.
When $d=O(\log n)$, we obtain a bound $2^{O(\log^2n)}$.
If in addition $\|L\|=O(1)$, then the  assumption (\ref{eq:degreealpha}) implies $\Delta=O(1)$, and then $T(d)=n^{O(1)}$.
\end{proof}

The following is then immediate.

\begin{coro}\label{coro:complexity}
Suppose $d=O(\log n)$. Then the complexity of the algorithm CountCOLOR is
$2^{O(\log^2n)}$. If in addition the size of the largest list $\|L\|$ is constant,
then CountCOLOR is a polynomial time algorithm.
\end{coro}

\section{Correlation decay}\label{section:CD}

The following is the key correlation decay result.

\begin{theorem}\label{theorem:correlationdecay}
Consider a triangle-free graph/list pair $(\G,\mbL)$ satisfying (\ref{eq:degreealpha}),(\ref{eq:beta}).
There exist constants $0<\epsilon<1$  which depend only on $\alpha$, such that for all nodes
$v$, colors $i\in L(v)$ and $d\geq 0$
\begin{align}\label{eq:correlationdecay}
\max_{i\in L(v)}
\Big|\log\pr_{\G,\mbL}(c(v)=i)-\log\Phi(\G,\mbL,v,i,d)\Big|\leq O(n^2(1-\epsilon)^d).
\end{align}
\end{theorem}
This theorem is our key tool for using the values of $\Phi$ for computing the marginals
$\pr_{\G,\mbL}(c(v)=i)$.
We first establish that this correlation decay result implies our main result, Theorem~\ref{theorem:MainResult}.

\begin{proof}[Proof of Theorem~\ref{theorem:MainResult}]
We consider an arbitrary instance $(\G,\mbL)$ with size $n$ and arbitrary $\delta>0$.
We may assume without the loss of generality that $n$ is at least a  large constant bigger than $C/\delta$,
for any universal constant $C$,
since we can simply extend the size of the instance by adding isolated nodes.
The proof uses a standard idea of approximating marginals $\pr_{\G,\mbL}(c(v)=i)$ and then using Proposition~\ref{prop:cavity}
for computing $Z(\G,\mbL)$. From Proposition~\ref{prop:cavity}, if the algorithm CountCOLOR produces in every stage
$k=1,2,\ldots,|V|-1$
a value $\hat p(v,i)$ which approximates $\pr_{\G_v,\mbL_k}(c(v_k)=i)$ with accuracy
\begin{align}\label{eq:accuracy1}
1-{\delta\over n}\leq {\hat p(v,i)\over \pr_{\G_v,\mbL_k}(c(v_k)=i)}\leq 1+{\delta\over n}
\end{align}
then the output $\hat Z$ of the algorithm satisfies
\begin{align*}
\Big(1-{\delta\over n}\Big)^{n}\leq
\Big(1-{\delta\over n}\Big)^{|V|}\leq{Z(\G,\mbL)\over \hat Z}\leq \Big(1+{\delta\over n}\Big)^{|V|}
\leq \Big(1+{\delta\over n}\Big)^{n}
\end{align*}
Since $|V|\leq n$ and $n$ is at least a large constant, we obtain an arbitrary accuracy of the approximation.
Thus it suffices to arrange for (\ref{eq:accuracy1}). We run the algorithm CountCOLOR with
$d=\lceil{4\log n\over \log {1\over 1-\epsilon}}\rceil$, where $\epsilon$ is the constant from Theorem~\ref{theorem:correlationdecay}.
This choice of $d$ gives $(1-\epsilon)^d\leq 1/n^4$.
Theorem~\ref{theorem:correlationdecay} with the given value of $d$ then implies
\begin{align*}
\Big|\log {\pr_{\hat \G,\hat \mbL}(c(v)=i)\over \hat p(v,i)}\Big|=
\Big|\log {\pr_{\hat \G,\hat \mbL}(c(v)=i)\over  \Phi(\hat \G,\hat \mbL,v,i,d)}\Big|\leq O(n^2){1\over n^4}=O({1\over n^2}).
\end{align*}
Thus
\begin{align*}
1-O({1\over n^2})\leq\exp\big(-O({1\over n^2})\big)
\leq {\hat p(v,i)\over \pr_{\G,\mbL}(c(v)=i)}\leq\exp\big(O({1\over n^2})\big)=1+O({1\over n^2})
\end{align*}
This gives us (\ref{eq:accuracy1}) for all $n>C/\delta$ where $C$ is the universal constant appearing in $O(\cdot)$.
This completes the analysis of the accuracy. The complexity part of the theorem follows directly from
Corollary~\ref{coro:complexity}.
\end{proof}

The rest  of the section is devoted to establishing this Theorem~\ref{theorem:correlationdecay}.
The basis of the proof is the recursion (\ref{eq:Aalpha}).
As before, let $v_1,\ldots,v_m$ be the neighbors of $v$ in $\G$, $m=\Delta(v)$.
Observe that (\ref{eq:correlationdecay}) holds trivially when $m=0$, since
both expression inside the absolute value become $1/|L(v)|$ and
the left-hand side becomes equal to zero.
Thus we assume that $m\geq 1$.
Denote by $m_k$ the degree of $v_k$ in the graph $\G_v$.
In order to ease the notations, we introduce
\begin{align*}
x_i&=\pr_{\G,\mbL}(c(v)=i),  ~~i\in L(v),\\
x_{i,k}&=\pr_{\G_v,\mbL_{k,i}}(c(v_k)=i),  ~~i\in L(v),1\leq k\leq m\\
x^*_i&=\Phi(\G,\mbL,v,i,d), ~~i\in L(v), \\
x^*_{i,k}&=\Phi(\G_v,\mbL_{k,i},v_k,i,d-1),  ~~i\in L(v)\cap L_{i,k}(v_k),1\leq k\leq m\\
\end{align*}
\ignore{
We have
\begin{align}
\sum_{i\in L(v)}x_i&=1, \label{eq:sumxj} \\
\sum_{i\in L(v)}x_i^*&\leq 1 \label{eq:sumxj'},
\end{align}
where (\ref{eq:sumxj'}) follows from Lemma~\ref{lemma:Phiproperties}.
}


\ignore{
and introduce
\begin{align*}
{\cal H}_{\G,\mbL}(\mbx,\mbx^*)=
\max_{j\in L(v)}\big|\log(x_j)-\log(x_j^*)\big|
\end{align*}
}

\begin{prop}\label{prop:H}
There exists a constant $\epsilon>0$ which depends on $\alpha$ only such that
\begin{align}\label{eq:recursion2}
{1\over m}\max_{i\in L(v)}\big|\log(x_i)-\log(x_i^*)\big|\leq (1-\epsilon)
\max_{j\in L(v),k:m_k>0}{1\over m_k}\big|\log(x_{j,k})-\log(x_{j,k}^*)\big|
\end{align}
\end{prop}

First we show how this result implies Theorem~\ref{theorem:correlationdecay}:
\begin{proof}[Proof of Theorem~\ref{theorem:correlationdecay}]
Applying this proposition $d$ times and using the fact that we are summing over $k:m_k>0$,
we obtain
\begin{align*}
{1\over m}\max_{i\in L(v)}\big[\log(x_i)-\log(x_i^*)\big]\leq M(1-\epsilon)^d,
\end{align*}
where
\begin{align*}
M=\max_{l,s}\Big|\log \pr_{\G_s,\mbL_s}(c(v)=l)-\log\Phi(\G_s,\mbL_s,v,l,0)\Big|
\end{align*}
and the maximum is over all graph/list pairs $\G_s,\mbL_s$ appearing during the computation of
$\Phi$ and over all colors $l$. Recall that if $l$ does not belong to the list associated with
node $v$ and list vector $\mbL_s$, then $\pr_{\G_s,\mbL_s}(c(v)=l)=\Phi(\G_s,\mbL_s,v,l,0)=0$
(the first is equal to zero by definition, the second by the way we set the values of $\Phi$). Otherwise we have
from Lemma~\ref{lemma:lessbeta} and part (\ref{eq:lowerboundPhi}) of Lemma~\ref{lemma:Phiproperties}
that absolute value of the difference is at most
\begin{align*}
\log q+\Delta \log(\beta/(\beta-1)).
\end{align*}
Since $m\leq \Delta\leq n$,  $\beta$ is a constant which only depends on $\alpha$, and $q\leq n$, then we obtain
$M=O(n)$ and $mM=O(n^2)$.
\end{proof}

Thus we focus on establishing Proposition~\ref{prop:H}.

\begin{proof}[Proof of Proposition~\ref{prop:H}]
Observe that for every $i\in L(v)\setminus L_{k,j}(v_k)$ we have
$x_{i,k}=x_{i,k}^*=0$. This is because the probability of node $v_k$ obtaining color $i$ is zero
when this color is not in its list. Similarly, the corresponding value of $\Phi$ is zero, since
we set it to zero for all colors not in the list.
For every $i\in L(v)$ introduce
\begin{align}\label{eq:Ak}
A_i\triangleq \prod_{1\leq k\leq m}(1-x_{i,k})
\end{align}
and
\begin{align}\label{eq:Asum}
A\triangleq\sum_{j\in L(v)}A_j
\end{align}
Introduce $A^*_i,A^*$ similarly. Applying Proposition~\ref{prop:Recursion1} we obtain
\begin{align}
&x_i={A_i\over A}, \label{eq:AoverA}\\
&x_i^*=\min\Big[{1\over 2(1+\epsilon_0)m},{1\over \beta},{A_i^*\over A^*}\Big].
\end{align}
Let
\begin{align*}
\tilde x_i^*={A_i^*\over A^*}.
\end{align*}
We claim that in order to establish (\ref{eq:recursion2})
it suffices to establish the bound
\begin{align*}
{1\over m}|\log x_i-\log \tilde x_i^*|\leq (1-\epsilon)
\max_{j\in L(v),k:m_k>0}{1\over m_k}\big|\log(x_{j,k})-\log(x_{j,k}^*)\big|
\end{align*}
Indeed, if $\tilde x_i^*\neq x_i^*$, then $x_i^*=\min[{1\over 2(1+\epsilon_0)m},{1\over \beta}]$.
On the other hand, by Lemmas~\ref{lemma:lessbeta},\ref{lemma:lessdelta} we have
$x_i\leq \min[{1\over 2(1+\epsilon_0)m},{1\over \beta}]$, implying $x_i\leq x_i^*\leq \tilde x_i^*$, and the
bound for $\tilde x_i^*$ implies the bound (\ref{eq:recursion2}).

We have
\begin{align}\label{eq:H1}
\max_{i\in L(v)}\big|\log(x_i)-\log(x_i^*)\big|&=\max_{i\in L(v)}\Big|\log A_i-\log A_i^*-\log A+\log A^*\Big|.
\end{align}
We introduce auxiliary variables $y_i=\log(x_i),y_{i,k}=\log(x_{i,k})$.
Similarly, let $y_i^*=\log(\tilde x_i^*),y_{i,k}^*=\log(x_{i,k}^*)$.
Define
$\mby=(y_{i,k}),\mby^*=(y_{i,k}^*)$. Observe that if $m_k=0$ then for every color $i$
$x_{i,k}=x_{i,k}^*$. This follows since both values are $1/|L_{i,k}|$ when $i\in L_{i,k}$ and
zero otherwise. This implies $y_{i,k}=y_{i,k}^*$.
Then we rewrite (\ref{eq:H1}) as
\begin{align}\label{eq:x}
\max_{i\in L(v)}\big|y_i-y_i^*\big|&=
\max_{i\in L(v)}\Big|\sum_{k:m_k>0}\log(1-\exp(y_{i,k}))-
\sum_{k:m_k>0}\log(1-\exp(y^*_{i,k}))\notag\\
&-\log\Big(\sum_{j\in L(v)}\prod_{1\leq k\leq m}(1-\exp(y_{j,k}))\Big)
+\log\Big(\sum_{j\in L(v)}\prod_{1\leq k\leq m}(1-\exp(y_{j,k}^*))\Big)
\Big|,
\end{align}
where  the sums $\sum_{1\leq k\leq m}$ were replaced by $\sum_{k:m_k>0}$ due to our observation
$y_{i,k}=y_{i,k}^*$ when $m_k=0$.

For every $i$ denote the expression inside the absolute value in the right-hand side of  equation (\ref{eq:x}) by
${\cal G}_i(\mby)$. That is  we treat $\mby^*$ as constant and $\mby$ as a variable.
It suffices to prove that for each $i$
\begin{align}\label{eq:Gi}
{\cal G}_i(\mby)\leq (1-\epsilon)
\max_{j\in L(v),k:m_k>0}{1\over m_k}\big|\log(x_{j,k})-\log(x_{j,k}^*)\big|
\end{align}
Observe that ${\cal G}_i(\mby^*)=0$.
Let $g_i(t)={\cal G}_i(\mby^*+t(\mby-\mby^*)), t\in [0,1]$. Then $g_i$ is a differentiable function interpolating between
$0$ and ${\cal G}_i(\mby)$. In particular, $g_i(1)={\cal G}_i(\mby)$.
Applying the Mean Value Theorem we obtain
\begin{align*}
|g_i(1)-g_i(0)|=|g_i(1)|&\leq \sup_{0\leq t\leq 1}|\dot g_i(t)| \\
&=\sup_{0\leq t\leq 1}\Big|\nabla {\cal G}_i(\mby^*+t(\mby-\mby^*))^T(\mby-\mby^*)\Big|\\
\end{align*}
where the supremum is over values of $t$. We use a short-hand notation
\begin{align*}
\Pi_j=\prod_{1\leq k\leq m}(1-\exp(y_{j,k}+t(y_{j,k}-y_{j,k}^*)))
\end{align*}
For each $t$ we have
\begin{align*}
&\nabla {\cal G}_i(\mby^*+t(\mby-\mby^*))(\mby-\mby^*)=
\sum_{k:m_k>0}{-\exp(y_{i,k}+t(y_{i,k}-y_{i,k}^*))\over 1-\exp(y_{i,k}+t(y_{i,k}-y_{i,k}^*))}
(y_{i,k}-y_{i,k}^*)\\
&+{\sum_{j\in L(v)}\sum_{1\leq k\leq m}{\exp(y_{j,k}+t(y_{j,k}-y_{j,k}^*))\over
1-\exp(y_{j,k}+t(y_{j,k}-y_{j,k}^*))}(y_{j,k}-y_{j,k}^*))
\Pi_j\over
\sum_{j\in L(v)}\Pi_j}.
\end{align*}
Again using the fact $y_{j,k}=y_{j,k}^*$ when $m_k=0$, we can replace the sum $\sum_{1\leq k\leq m}$ by
$\sum_{k:m_k>0}$ in the expression above.
For each $j$ we have from convexity of $\exp$
\begin{align*}
\exp(y_{j,k}+t(y_{j,k}-y_{j,k}^*))&\leq (1-t)\exp(y_{j,k}^*)+t\exp(y_{j,k})\\
&=(1-t)x_{j,k}+t x_{j,k}^* \\
&\leq {1\over 2(1+\epsilon_0)m_k}.
\end{align*}
where the last inequality follows from Lemma~\ref{lemma:lessdelta} and
part (\ref{eq:PhiProperties1}) of Lemma~\ref{lemma:Phiproperties}.
This bound is useful for terms with $m_k>0$ (for this reason we only kept these terms in the sum $\sum_{k:m_k>0}$).
Similarly using Lemma~\ref{lemma:lessbeta} and again part (\ref{eq:PhiProperties1}) of Lemma~\ref{lemma:Phiproperties}
we obtain
\begin{align*}
{1\over 1-\exp(y_{j,k}+t(y_{j,k}-y_{j,k}^*))}&\leq {1\over 1-(1-t)\exp(y_{j,k}^*)-t\exp(y_{j,k})}\\
&={1\over 1-(1-t)x_{j,k}-t x_{j,k}^*} \\
&\leq {1\over 1-{1\over \beta}}.
\end{align*}
We obtain
\begin{align*}
\sup_{0\leq t\leq 1}\Big|\nabla {\cal G}_i(\mby^*+t(\mby-\mby^*))(\mby-\mby^*)\Big|
&\leq
\sum_{k:m_k>0}{1\over(1-{1\over \beta})2(1+\epsilon_0)m_k}
|y_{i,k}-y_{i,k}^*|\\
&+{\sum_{j\in L(v)}\sum_{k:m_k>0}(1-{1\over \beta})^{-1}2^{-1}(1+\epsilon_0)^{-1}m_k^{-1}|y_{j,k}-y_{j,k}^*|
\Pi_j\over
\sum_{j\in L(v)}\Pi_j} \\
&\leq {m\over(1-{1\over \beta})2(1+\epsilon_0)}\max_{j\in L(v),k:m_k>0}{|y_{j,k}-y_{j,k}^*|\over m_k}\\
&+{m\over(1-{1\over \beta})2(1+\epsilon_0)}\max_{j\in L(v),k:m_k>0}{|y_{j,k}-y_{j,k}^*|\over m_k}\\
&= {m\over(1-{1\over \beta})(1+\epsilon_0)}\max_{j\in L(v),k:m_k>0}{|y_{j,k}-y_{j,k}^*|\over m_k}.
\end{align*}
Combining with (\ref{eq:x}) we conclude
\begin{align*}
\max_{i\in L(v)}{|y_i-y_i^*|\over m}\leq
{1\over (1-{1\over \beta})(1+\epsilon_0)}\max_{j\in L(v),k:m_k>0}{|y_{j,k}-y_{j,k}^*|\over m_k}.
\end{align*}

We now select a sufficiently large constant $\beta=\beta(\epsilon_0)$ such that
\begin{align*}
1-\epsilon\triangleq {1\over (1-{1\over \beta})(1+\epsilon_0)}<1.
\end{align*}
This completes the proof of Proposition~\ref{prop:H}.
\end{proof}

\section{Extensions:  Markov random field and partition function}
\label{section:MRF}

\subsection{Model and the preliminary results}
The main conceptual point of this paper, namely construction a recursion of the form (\ref{eq:Aalpha}),
construction of a corresponding  computation tree, establishing correlation decay property and application
to a counting problem,
can be
extended to an arbitrary model of random constraint satisfaction problems with multiple values.
In this section we provide details using a very general framework of Markov random fields (MRF),
also known as \emph{graphical model} \cite{JordanGraphicalModels},\cite{JordanWainwright}.
We show that generalizing (\ref{eq:Aalpha}) is straightforward. It is establishing
the decay of correlation which presents the main technical difficulty. We provide a simple and general
sufficient condition and then illustrate the approach on specific statistical physics problem, namely
$q$-state Potts model. Here we restrict ourselves for simplicity to MRF defined on simple graphs. Extensions to
multi-graphs are possible as well.

A Markov random field (MRF) is given as a graph $\G$ with node set $V=\{v_1,\ldots,v_{|V|}\}$,
edge set $E$, an alphabet $\X$
and set of functions $\phi_v:\X:\rightarrow \R_+, ~v\in V,~ f_{v,u}:\X\rightarrow \R_+, (v,u)\in E$.
Consider a probability measure on $\X^{|V|}$ defined by
\begin{align*}
\pr(\mbX=\mbx)={\prod_{v\in V}\phi_v(x_v)\prod_{(v,u)\in E}f_{u,v}(x_v,x_u)\over Z},
\end{align*}
for every $\mbx=(x_v)\in \X^{|V|}$, where $Z=\sum_{\mbx}\prod_{v\in V}\phi_v(x_v)\prod_{(v,u)\in E}f_{v,u}(x_v,x_u)$
is the normalizing constant called the partition function.
Here  $\mbX=(X_v)$ is the random vector selected according to this probability measure.
In the case $Z=0$, the MRF is not defined.
From now on assume that \\
$\prod_{v\in V}\phi_v(x_v)\prod_{(v,u)\in E}f_{v,u}(x_v,x_u)>0$ for at least
one $\mbx=(x_v)\in \X^{|V|}$.

Let us see that the problem of list-coloring can be cast as a Markov random field, where $\pr(\cdot)$ corresponds
to the uniform probability distribution on the set of valid colorings. Given an instance of a list-coloring
problem $(\G,\mbL)$ with a universe of colors $\{1,\ldots,q\}$,
we set $\X=\{1,\ldots,q\}$, $\phi_v(i)=1\{i\in L(v)\}$ for all node/color pairs $v,i$,
and $f_{v,u}(i,j)=1\{i\ne j\}$, where $1\{\cdot\}$ is the indicator function. It is not hard to see that
$\pr(\mbx)=1/Z$ if $\mbx$ corresponds to a valid coloring and $=0$ otherwise, and $Z=Z(\G,\mbL)$ is the total
number of valid list-colorings. Thus this MRF corresponds to the uniform distribution on the set of proper colorings.

An instance of a MRF is denoted by $\M=(\G,\X,\phi,f)$, with $\phi=(\phi_v),f=(f_{v,u})$.
We will write $\pr_{\M}$ and $Z_{\M}$ for the corresponding probability measure and the partition function, respectively,
in order to emphasize the dependence on the particular instance of the MRF. Computation of $Z_{\M}$ is the principle
goal of this section.
As in the case of list-coloring model, denote by $Z_{\M}[\chi]$ the sum of the terms in the partition function
which satisfy some condition $\chi$.

Observe that if $\phi_v,f_{v,u}>0$ for all nodes and edges than $\pr_{\M}(\mbX=\mbx)>0$ for every  $\mbx=(x_v), v\in V$.
Moreover, if $f_{v,u}=c$ for all edges for some constant $c$, then
we obtain a product form solution
\begin{align*}
\pr_{\M}(\mbX=\mbx)=\prod_v{\phi_v(x_v)\over \sum_{y\in\X}\phi_v(y)}.
\end{align*}
Thus we might expect
the correlation decay to take place when the values of $f_{v,u}$ are \emph{close} to each other. This is the regime within
which we will establish our results. Let $\phi_{\min}=\min_{v,x}\phi_v(x),\phi_{\max}=\max_{v,x}\phi_v(x)$
and $c_{\phi}=\phi_{\max}/\phi_{\min}$.
Also let
$f_{\min}=\min_{(v,u)\in E, x,y\in\X}f_{v,u}(x,y),
f_{\max}=\max_{(v,u)\in E, x,y\in\X}f_{v,u}(x,y)$ and let $c_f=f_{\max}/f_{\min}$.
From now on we assume that the following conditions hold
\begin{align}
\phi_v(x)&>0,~~\forall v\in V, x\in\X \label{eq:phinonnegative}\\
f_{v,u}(x,y)&>0, ~~\forall (v,u)\in E, ~x,y\in \X  \label{eq:fnonnegative}
\end{align}
These conditions in particular ensure that $c_f<\infty$. The following condition will be used in lieu  of (\ref{eq:degreealpha})
\begin{align}
\gamma\triangleq (c_f^\Delta-c_f^{-\Delta})\Delta |X|^\Delta< 1. \label{eq:conditionMRF}
\end{align}
The size of an instance
$\M$ is
\begin{align*}
n=\max\Big(|V|,|E|,|\X|,|\log \phi_{\max}|,|\log \phi^{-1}_{\min}|,|\log f_{\max}|,|\log f^{-1}_{\min}|\Big).
\end{align*}

We now state the main result of this section.
\begin{theorem}\label{theorem:MainResultMRF}
There exist a deterministic algorithm which provides an FPTAS for computing $Z_\M$ for an arbitrary
MRF instance  $\M$ satisfying (\ref{eq:phinonnegative}),(\ref{eq:fnonnegative}),(\ref{eq:conditionMRF}),
whenever   $|\X|$ and $\Delta$ are
constants.
\end{theorem}

Our first task is obtaining a generalization of the cavity recursion given by Proposition~\ref{prop:cavity}.
Given a MRF $\M=(\G,\X,\phi,f)$, an arbitrary node $v$ and an arbitrary element $x^*\in\X$ we consider a new MRF instance
$\mathcal{T}_{v,x^*}[\M]=(\tilde\G,\tilde \X,\tilde \phi,\tilde f)$ defined as follows.
The graph $\tilde\G$ is the subgraph of $\G$
induced by all nodes other than $v$. $\tilde \X=\X$ and $\tilde f=f$. $\tilde \phi$ is defined as follows.
For every $u$ which is a neighbor of $v$, $\tilde \phi_u(x)=\phi_v^{1\over \Delta(v)}f_{v,u}(x^*,x)$,
where $\Delta(v)$ is the degree of $v$ in $\G$. For all the remaining nodes $u$ we set $\tilde \phi_u=\phi_u$.

Given a MRF $\M=(\G,\X,\phi,f)$ let $v_1,\ldots,v_{|V|}$ be an arbitrary enumeration of nodes.
Consider an arbitrary $\mbx^*=(x_{v_1}^*,\ldots,x_{v_{|V|}}^*)$ such that $\pr(\mbx^*)>0$.
Define $\M_0=\M$ and $\M_k=\mathcal{T}_{v_k,x^*_{v_k}}[\M_{k-1}], ~k=1,2,\ldots,|V|$,
where $\M_{|V|}$ is an empty MRF and its partition function is set by default to  unity.

\begin{prop}\label{prop:cavityMRF}
 The following identity holds.
\begin{align*}
Z_{\M}=\prod_{1\leq k\leq |V|}\pr_{\M_{k-1}}^{-1}(X_{v_k}=x_{v_k}^*).
\end{align*}
\end{prop}

\begin{proof}
We have
\begin{align*}
\pr_{\M}(X_{v_1}=x^*_{v_1})&={\sum_{\mbx\in \X^{|V|}:x_{v_1}=x^*_{v_1}}\phi_{v_1}(x^*_{v_1})
\prod_{u:(v_1,u)\in E}f_{v_1,u}(x^*_{v_1},x_u)
\prod_{u\ne v_1}\phi(x_u)\prod_{(v,u)\in E,~v,u\ne v_1}f_{v,u}(x_v,x_u)\over Z_{\M}} \\
&={\prod_{u\ne v_1}\phi_u^1(x_u)\prod_{(v,u)\in E,~v,u\ne v_1}f_{v,u}(x_v,x_u)\over Z_{\M}} \\
&={Z_{\M_1}\over Z_{\M}},
\end{align*}
where the second equality follows since $\phi(x^*_{v_1})\prod_{u:(v_1,u)\in E}f(x^*_{v_1},x_u)=\prod_{u}\phi^1(x_u)$
and the second product is over neighbors $u$ of $v_1$ in $\G$.
Iterating further for $k\geq 2$ we obtain the result.
\end{proof}

The identity in Proposition~\ref{prop:cavityMRF} provides an important representation  of the partition function
in terms of marginal probabilities. Thus, if we compute (approximately) these marginal probabilities, we can use
them to obtain the value of the underlying partition function.

\subsection{Basic recursion and the algorithm}
Our next task is constructing a generalization of $(\G_v,\mbL_{k,i})$ and extending Proposition~\ref{prop:Recursion1}
to MRF. Given a MRF $\M=(\G,\X,\phi,f)$ and a node $v$ let $\M_v$ denote the MRF instance obtained naturally
by removing node $v$. Namely, we keep $\phi_u$ and $f_{u,w}$ intact for all the nodes $u\ne v$ and edges $(u,w), u,w\ne v$.
Also, given a MRF $\M=(\G,\X,\phi,f)$,
a set of nodes $v_1,\ldots,v_r\subset V$ and a set of
elements $x_1,\ldots,x_r\in \X$ we construct a MRF
denoted by $\M[v_1,x_1;\ldots; v_r,x_r]=(\tilde \G,\tilde\phi,\tilde f)$ as follows.
The corresponding graph $\tilde \G$ is the subgraph induced by nodes $V\setminus \{v_1,\ldots,v_r\}$.
For every node $u\in \tilde \G$ which has at least one neighbor among $v_1,\ldots,v_r$
we set $\tilde\phi_u(x)=\prod_{i}f_{v_i,u}(x_i,x)\phi_u(x)$,
where the product is over $i=1,2\ldots,r$ such that $(v_i,u)$ is an edge in $\G$. For all the remaining
$u$ we set $\tilde\phi_u=\phi_u$. We also set $\tilde f=f$. The interpretation for $\M[v_1,x_1;\ldots; v_r,x_r]$
comes from the following simple fact.

\begin{lemma}\label{lemma:CondM}
For every event $\mathcal{E}$ corresponding to the probability measure $\pr_\M$, the following holds
\begin{align*}
\pr_\M(\mathcal{E}| \wedge_{k\le r}X_{v_k}=x_k)=\pr_{\M[v_1,x_1;\ldots; v_r,x_r]}(\mathcal{E}).
\end{align*}
\end{lemma}

\begin{proof}
The proof is obtained immediately by summing over all of the elementary events $\mbx\in\mathcal{E}$ and
observing that the terms $\phi_{v_k}(x_k)$ cancel in the ratio
$\pr_M(\mathcal{E}\wedge \wedge_{k\le r}X_{v_k}=x_k)/\pr_M(\wedge_{k\le r}X_{v_k}=x_k)$.
\end{proof}
Observe that the value of $c_f$ corresponding to the MRF $\M[v_1,x_1;\ldots; v_r,x_r]$
is the same of $\M$. Thus,
should $\M$ satisfy conditions
(\ref{eq:phinonnegative}),(\ref{eq:fnonnegative}),(\ref{eq:conditionMRF}),
so does the instance $\M[v_1,x_1;\ldots; v_r,x_r]$. Moreover, $c_{\tilde\phi}$ defined for this MRF
satisfies
\begin{align}\label{eq:ctildephi}
c_{\tilde\phi}\le c_{\phi}c_f^{\Delta}
\end{align}

Now we obtain a recursion which serves as a basis for our correlation decay analysis and construction of an algorithm.
\begin{prop}\label{prop:Recursion1MRF}
For every  node $v$ and its neighbors $v_1,\ldots,v_m$, the
following identity holds for every $x_0\in \X$:
\begin{align}
\pr_{\M}(X_v=x_0)= {\phi_v(x_0)\sum_{x_1,\ldots,x_m\in \X
}\prod_{k=1}^mf_{v,v_k}(x_0,x_k)\pr_{\M[v_1,x_1;\ldots;v_{k-1},x_{k-1}]}(X_{v_k}=x_k)
\over \sum_{x\in\X}\phi_v(x)\sum_{x_1,\ldots,x_m\in \X }
\prod_{k=1}^mf_{v,v_k}(x,x_k)\pr_{\M[v_1,x_1;\ldots;v_{k-1},x_{k-1}]}(X_{v_k}=x_k)},
\label{eq:Recursion1MRF}
\end{align}
where the sum $\sum_{x_1,\ldots,x_m\in \X }=1$ when $m=0$.
\end{prop}
\begin{proof}
The case $m=0$ is immediate. Assume $m>0$.
For every $x_0\in\X$ we have the following identity
\begin{align*}
\pr_{\M}(X_v=x_0)= {\phi_v(x_0)\sum_{x_1,\ldots,x_m\in \X
}Z_{\M_v}[X_{v_1}=x_1,\ldots,X_{v_m}=x_m]\prod_{k=1}^mf_{v,v_k}(x_0,x_k)
\over \sum_{x\in\X}\phi_v(x)\sum_{x_1,\ldots,x_m\in \X
}Z_{\M_v}[X_{v_1}=x_1,\ldots,X_{v_m}=x_m]\prod_{k=1}^mf_{v,v_k}(x,x_k)}
\end{align*}
We divide both parts by $Z_{\M_v}$ and write
\begin{align*}
{Z_{\M_v}[X_{v_1}=x_1,\ldots,X_{v_m}=x_m]\over Z_{\M_v}}=
\prod_{k=1}^m{Z_{\M_v}[X_{v_1}=x_1,\ldots,X_{v_k}=x_k]\over Z_{\M_v}[X_{v_1}=x_1,\ldots,X_{v_{k-1}}=x_{k-1}]},
\end{align*}
where the  term corresponding to $k=0$ is identified with $Z_{\M_v}$.
Applying Lemma~\ref{lemma:CondM}, we recognize the $k$-th term in this product as
$\pr_{\M[v_1,x_1;\ldots;v_{k-1},x_{k-1}]}(X_{v_k}=x_k)$ (note that the terms $\phi_{v_j}(x_j), j\le k-1$ cancel out).

\end{proof}

Proposition~\ref{prop:Recursion1MRF} also allows us to obtain
upper and lower bounds on the marginal probabilities:
\begin{lemma}\label{lemma:boundphi}
For every node $v$ and $x_0\in \X$
\begin{align*}
c_f^{-\Delta}{\phi_v(x_0)\over \sum_x \phi_v(x)}\le \pr_{\M}(X_{v}=x_0)
\le c_f^{\Delta}{\phi_v(x_0)\over \sum_x \phi_v(x)}.
\end{align*}
\end{lemma}

\begin{proof}
The proof follows from Proposition~\ref{prop:Recursion1MRF}. We have for every $x\in \X$, node $v$ and
its neighbors $v_1,\ldots,v_m$ that
$f_{\min}^m\le \prod_{k=1}^mf_{v,v_k}(x,x_k)\le c_f^m f_{\min}^m$. Applying this bound to the numerator
of (\ref{eq:Recursion1MRF}) for $x=x_0$ we obtain the required upper bound. Applying the same to the denominator,
we obtain the required lower bound.
\end{proof}

We now provide sufficient conditions under which the construction of a computation tree for computing approximately
marginal probabilities $\pr_{\M}(X_v=x)$ as well as the partition function $Z_{\M}$ can be performed in polynomial time.

Similarly to the problem of coloring, we introduce $\Phi(\cdot)$ -- a surrogate for computing the marginal probabilities
$\pr_{\M}(\cdot)$. Consider a function $\Phi_{\M}(v,x,d)$ defined recursively for an arbitrary instance of a MRF $\M=(\G,\X,\phi,f)$,
arbitrary node $v$, element $x\in \X$ and a non-negative integer $d$ as follows.

\begin{itemize}
\item We set $\Phi_{\M}(v,x,0)=1$. As in the case of coloring,
it turns out that the initialization values are not particularly important, due to the decay of correlations.

\item For every node $v$ with neighbors $v_1,\ldots,v_m$, every $x_0\in\X$ and $d\ge 1$
\begin{align}
\Phi_{\M}(v,x_0,d)= {\phi_v(x_0)\sum_{x_1,\ldots,x_m\in \X
}\prod_{k=1}^m\Phi_{\M[v_1,x_1;\ldots;v_{k-1},x_{k-1}]}(v_k,x_k,d-1)f_{v,v_k}(x_0,x_k)
\over \sum_{x\in\X}\phi_v(x)\sum_{x_1,\ldots,x_m\in \X
}\prod_{k=1}^m\Phi_{\M[v_1,x_1;\ldots;v_{k-1},x_{k-1}]}(v_k,x_k,d-1)
f_{v,v_k}(x,x_k)}, \label{eq:RecursionPhi}
\end{align}
where the sum $\sum_{x_1,\ldots,x_m\in \X }=1$ when $m=0$.
\end{itemize}
Assumptions (\ref{eq:phinonnegative}) and (\ref{eq:fnonnegative}) guarantee that $\Phi>0$.

We now describe our algorithm for approximately computing $Z_{\M}$.
The algorithm is parametrized by  $d$. It is based on computing recursively the values of $\Phi_{\M}$.

\vspace{.1in}

\textbf{Algorithm ComputeZ}
\vspace{.1in}

{\tt
INPUT: A MRF instance $\M=(\G,\X,\phi,f)$ and a positive integer $d$.

BEGIN

Set $\hat Z=1, \hat \M=\M$.

While $\hat \G\neq \emptyset$, fix an arbitrary node $v\in \hat \G$ and element $x\in \X$. Compute
$\Phi_{\hat\M}(v,x,d)$.

Set $\hat Z=\Phi^{-1}_{\hat\M}(v,x,d)\hat Z$.

Set $\hat\M=\mathcal{T}_{v,x}[\hat\M]$, where the operator $\mathcal{T}$
was defined before Proposition~\ref{prop:cavityMRF}.

END

OUTPUT: $\hat Z$.}

\vspace{.2in}

\subsection{Complexity}
We begin by analyzing the complexity of computing function $\Phi$.

\begin{prop}
For every $v\in V, x\in\X$, the function $\Phi_\M(v,x,d)$ can be computed in time $O(2^{d\Delta \log|\X|}n^2)$.
In particular when $d=O(\log n)$, and $|\X|,\Delta=O(1)$, the computation is polynomial in $n$.
\end{prop}
We note that the dependence on $\Delta$ is not as nice as in the case of the list-coloring problem, as
it appears as $\Delta$ not $\log\Delta$ in the exponent. Thus we can no longer claim that the computation time
is $2^{O(\log^2 n)}$ in this case.

\begin{proof}
Let $T(d)$ denote the complexity of computing function $\Phi(\cdot,d)$. Clearly,
$T(0)=O(1)$.
We now express $T(d)$ in terms of $T(d-1)$. Given a node $v$,
in order to compute $\Phi_{\M}(v,x,d)$ we identify the neighbors $v_1,\ldots,v_m$ of $v$.
For every sequence $x_1,\ldots,x_m\in \X$ and every $k=1,2,\ldots,m$
we compute $\Phi_{\M[v_1,x_1;\ldots;v_{k-1},x_{k-1}]}(v_k,x_k,d-1)$. The computation of each such quantity
is $T(d-1)$. We use the obtained values
to compute $\Phi_{\M}(v,x,d)$ via (\ref{eq:RecursionPhi}). We also need $O(n^2)$ time to "take care" of
multiplying by $f_{v,v_k}$ and by $\phi_v$.
The overall computation effort then satisfies
\begin{align*}
T(d)&=O(|\X|^\Delta T(d-1)+n^2).  
\end{align*}
Iterating over $d$ we obtain $T(d)=O(|\X|^{d\Delta}n^2)=O(2^{d\Delta\log|\X|}n^2)$, and the first part is established.
When $d=O(\log n)$ and $\Delta,|\X|$ are constants, we obtain $T(d)=n^{O(1)}$.
\end{proof}

The following is then immediate.

\begin{coro}\label{coro:complexityMRF}
Suppose $d=O(\log n)$ and $|\X|,\Delta=O(1)$. Then  ComputeZ is a
polynomial time algorithm.
\end{coro}

\subsection{Correlation decay analysis}
We now establish a correlation decay result which is  a key to proving our main result, Theorem~\ref{theorem:MainResultMRF}.

\begin{theorem}\label{theorem:correlationdecayMRF}
Given an arbitrary MRF satisfying conditions (\ref{eq:phinonnegative}),(\ref{eq:fnonnegative}),(\ref{eq:conditionMRF}),
the following holds for every node $v$ and $d\ge 1$
\begin{align}
\max_{x\in \X}\Big|&\log\pr_{\M}(X_v=x)-\log\Phi_{\M}(v,x,d)\Big| \notag\\
&\le (1-\gamma) \max_{1\le k\le m, y\in \X}
\Big|\pr_{\M[v_1,x_1;\ldots;v_{k-1},x_{k-1}]}(X_{v_k}=y)-\Phi_{\M[v_1,x_1;\ldots;v_{k-1},x_{k-1}]}(v_k,y,d-1)\Big|
\label{eq:contractionMRF}
\end{align}
where $v_1,\ldots,v_m$ are the neighbors of $v$.
\end{theorem}
We first show how this theorem implies our main algorithmic result.

\begin{proof}[Proof of Theorem~\ref{theorem:MainResultMRF}]
We claim that ComputeZ provides FPTAS for computing partition function $Z_{\M}$ when $d=O(\log n)$
under the setting of Theorem~\ref{theorem:MainResultMRF}.
We have already established in Corollary~\ref{coro:complexityMRF} that the algorithm is polynomial time.

Consider any MRF instance $\tilde \M$ obtained during the computation of $\Phi_\M(\cdot)$ as a part of performing
algorithm ComputeZ. Applying (\ref{eq:ctildephi}) and Lemma~\ref{lemma:boundphi} we obtain that for every
node $v$ in $\tilde M$ and every $x\in\X$
\begin{align*}
\pr_{\tilde\M}(X_v=x)\ge c_f^{-\Delta}{1\over |\X|}c_{\tilde \phi}^{-1}\ge
c_f^{-\Delta}{1\over |\X|}c_f^{-\Delta d}c_{\phi}^{-1}=c_f^{-\Delta(d+1)}{1\over |\X|}c_{\phi}^{-1}.
\end{align*}
Then applying the result of Theorem~\ref{theorem:correlationdecayMRF} $d$ times and
recalling $\Phi_{\tilde M}(v,x,0)=1$, we obtain for $d=O(\log n)$
\begin{align*}
\Big|\log\pr_{\M}(X_v=x)-\log\Phi_{\M}(v,x,d)\Big|&\le (1-\gamma)^d\Big(\Delta(d+1)\log c_f+\log|\X|+\log c_{\phi}\Big)\\
&=(1-\gamma)^d O(dn\log n)\\
&={1\over n^{O(\log {1\over 1-\gamma})}}O(n\log^2 n)\\
&={1\over n^{O(1)}}
\end{align*}
where the last step is obtained by selecting $d=C\log n$ for sufficiently large constant $C$.
\begin{align*}
\Big|{\pr_{\M}(X_v=x)\over \Phi_{\M}(v,x,d)}-1\Big|\le \exp(n^{-\Omega(1)})-1={1\over n^{\Omega(1)}}.
\end{align*}
We conclude that $\Phi_{\M}(v,x,d)$ provides an approximation of marginal probability $\pr_{\M}(X_v=x)$
with an inverse polynomial error. The remainder of the proof is the same as for Theorem~\ref{section:MainResult}.
\end{proof}

\begin{proof}[Proof of Theorem~\ref{theorem:correlationdecayMRF}]
Fix a node $v$ and an element $x_0\in\X$. Let $v_1,\ldots,v_m$ be neighbors of $v$.
When $m=0$ we the left-hand side of (\ref{eq:contractionMRF}) is zero. Thus assume $m>0$.
In order to ease the exposition we introduce some notations.
Set $z=\log\pr_{\M}(X_v=x_0)$, $z_{x_1,\ldots,x_k}=\log\pr_{\M[v_1,x_1;\ldots;v_{k-1},x_{k-1}]}(X_{v_k}=x_k)$.
Similarly $\tilde z=\log\Phi_{\M}(v,x_0,d), \tilde z_{x_1,\ldots,x_k}=\log\Phi_{\M[v_1,x_1;\ldots;v_{k-1},x_{k-1}]}(v_k,x_k,d-1)$.
Also let $\mbz$ denote the vector $(z_{x_1,\ldots,x_k}),1\le k\le m, x_1,\ldots,x_m\in \X$ and
$\tilde \mbz$ denote the vector $(\tilde z_{x_1,\ldots,x_k}),1\le k\le m, x_1,\ldots,x_m\in \X$.
Both vectors have dimension $\sum_{1\le k\le m}|\X|^m$.
Then we can rewrite (\ref{eq:Recursion1MRF}) as
\begin{align}
z=\log{\phi_v(x_0)\sum_{x_1,\ldots,x_m\in \X
}\prod_{k=1}^mf_{v,v_k}(x_0,x_k)\exp(z_{x_1,\ldots,x_k}) \over
\sum_{x\in\X}\phi_v(x)\sum_{x_1,\ldots,x_m\in \X }
\prod_{k=1}^mf_{v,v_k}(x,x_k)\exp(z_{x_1,\ldots,x_k})},
\label{eq:Recursion1MRFz}
\end{align}
and rewrite (\ref{eq:RecursionPhi}) as
\begin{align}
\tilde z=\log{\phi_v(x_0)\sum_{x_1,\ldots,x_m\in \X
}\prod_{k=1}^mf_{v,v_k}(x_0,x_k)\exp(\tilde z_{x_1,\ldots,x_k})
\over \sum_{x\in\X}\phi_v(x)\sum_{x_1,\ldots,x_m\in \X }
\prod_{k=1}^mf_{v,v_k}(x,x_k)\exp(\tilde z_{x_1,\ldots,x_k})},
\label{eq:Recursion1MRFtildez}
\end{align}
Introduce a function ${\cal G}$ defined on a vector $\mbw=(w_{x_1,\ldots,x_k}), 1\le k\le m, x_1,\ldots,x_m\in \X$
with the same dimension $\sum_{1\le k\le m}|\X|^m$
as follows:
\begin{align}
{\cal G}(\mbw)=\log{\phi_v(x_0)\sum_{x_1,\ldots,x_m\in \X
}\prod_{k=1}^mf_{v,v_k}(x_0,x_k) \exp(w_{x_1,\ldots,x_k}) \over
\sum_{x\in\X}\phi_v(x)\sum_{x_1,\ldots,x_m\in \X }
\prod_{k=1}^mf_{v,v_k}(x,x_k)\exp(w_{x_1,\ldots,x_k})},
\label{eq:calG}
\end{align}
which we rewrite as
\begin{align*}
\log\phi_v(x_0)+\log{\cal G}_1(\mbw)-\log{\cal G}_2(\mbw)
\end{align*}
where the definition of ${\cal G}_1$ and ${\cal G}_2$ is immediate.
\begin{align*}
{\cal G}_1(\mbw)&=\phi_v(x_0)\sum_{x_1,\ldots,x_m\in \X
}\prod_{k=1}^mf_{v,v_k}(x_0,x_k)
\exp(w_{x_1,\ldots,x_k}) \\
{\cal G}_2(\mbw)&=\sum_{x\in\X}\phi_v(x)\sum_{x_1,\ldots,x_m\in \X }
\prod_{k=1}^mf_{v,v_k}(x,x_k)\exp(w_{x_1,\ldots,x_k})
\end{align*}
We have $z-\tilde z={\cal G}(\mbz)-{\cal G}(\tilde \mbz)$. Thus establishing (\ref{eq:contractionMRF})
reduces to showing
\begin{align*}
|{\cal G}(\mbz)-{\cal G}(\tilde \mbz)|\leq (1-\gamma)\|\mbz-\tilde\mbz\|_{L_\infty}.
\end{align*}
Applying Mean Value Theorem, there exists $t\in [0,1]$ such that
\begin{align*}
z-\tilde z=\nabla {\cal G}(t\mbz+(1-t)\tilde \mbz)^T(\mbz-\tilde\mbz)
\end{align*}
further implying
\begin{align*}
|z-\tilde z|\le \|\nabla {\cal G}(t\mbz+(1-t)\tilde \mbz)\|_{L_1}\|\mbz-\tilde\mbz)\|_{L_\infty}.
\end{align*}
It then suffices to establish
\begin{align*}
\|\nabla {\cal G}(t\mbz+(1-t)\tilde \mbz)\|_{L_1}\le 1-\gamma.
\end{align*}
In the following lemma we show that this bound holds for an arbitrary input vector $\mbw$ and thus complete
the proof of Theorem~\ref{theorem:correlationdecayMRF}.
\end{proof}

\begin{lemma}
For every vector $\mbw$
\begin{align*}
\|\nabla {\cal G}(\mbw)\|_{L_1}\le 1-\gamma.
\end{align*}
\end{lemma}

\begin{proof}
Fix an arbitrary sequence $x^0_1,\ldots,x^0_{k_0}\in \X$ and the corresponding variable $w_{x^0_1,\ldots,x^0_{k_0}}$.
We have
\begin{align*}
{\partial {\cal G}\over \partial w_{x^0_1,\ldots,x^0_{k_0}}}=
{\cal G}_1^{-1}{\partial {\cal G}_1\over \partial w_{x^0_1,\ldots,x^0_{k_0}}}-
{\cal G}_2^{-1}{\partial {\cal G}_2\over \partial w_{x^0_1,\ldots,x^0_{k_0}}}
\end{align*}
We have
\begin{align*}
&{\cal G}_1^{-1}{\partial {\cal G}_1\over \partial w_{x^0_1,\ldots,x^0_{k_0}}}= \\
&{\Big(\prod_{k=1}^{k_0}f_{v,v_k}(x_0,x^0_k)\exp(w_{x^0_1,\ldots,x^0_k})\Big)
\sum_{x_{k+1},\ldots,x_m\in \X}\prod_{k=k_0+1}^mf_{v,v_k}(x_0,x_k)
\exp(w_{x^0_1,\ldots,x^0_{k_0},x_{k_0+1}\ldots,x_k})\over
\sum_{x_1,\ldots,x_m\in \X}\prod_{k=1}^mf_{v,v_k}(x_0,x_k)
\exp(w_{x_1,\ldots,x_k})}
\end{align*}
Using $f_{\min}\le f_{v,v_k}(x_0,x^0_k)\leq c_f f_{\min}$, we obtain
\begin{align*}
{\cal G}_1^{-1}{\partial {\cal G}_1\over \partial w_{x^0_1,\ldots,x^0_{k_0}}}
&\le c_f^m{\Big(\prod_{k=1}^{k_0}\exp(w_{x^0_1,\ldots,x^0_k})\Big)
\sum_{x_{k+1},\ldots,x_m\in \X}\prod_{k=k_0+1}^m
\exp(w_{x_1,\ldots,x_k})\over
\sum_{x_1,\ldots,x_m\in \X}\prod_{k=1}^m\exp(w_{x_1,\ldots,x_k})} \\
&\le c_f^{m} \\
&\le c_f^\Delta.
\end{align*}
Similarly, we obtain
\begin{align*}
{\cal G}_1^{-1}{\partial {\cal G}_1\over \partial w_{x_1,\ldots,x_k}}\ge  c_f^{-\Delta}.
\end{align*}
Using again $f_{\min}\le f_{v,u}(x,y)\leq c_f f_{\min}$ we also obtain
\begin{align*}
&{\cal G}_2^{-1}{\partial {\cal G}_2\over \partial w_{x^0_1,\ldots,x^0_{k_0}}}=\\
&={\sum_{x\in\X}\phi_v(x)\Big(\prod_{k=1}^{k_0}f_{v,v_k}(x,x^0_k)\exp(w_{x^0_1,\ldots,x^0_k})\Big)
\sum_{x_{k+1},\ldots,x_m\in \X}\prod_{k=k_0+1}^mf_{v,v_k}(x_0,x_k)
\exp(w_{x^0_1,\ldots,x^0_{k_0},x_{k_0+1}\ldots,x_k})\over
\sum_{x\in\X}\phi_v(x)\sum_{x_1,\ldots,x_m\in \X}\prod_{k=1}^mf_{v,v_k}(x,x_k)
\exp(w_{x_1,\ldots,x_k})}\\
&\le c_f^m
{\Big(\sum_{x\in\X}\phi_v(x)\Big)\Big(\prod_{k=1}^{k_0}\exp(w_{x^0_1,\ldots,x^0_k})\Big)
\sum_{x_{k+1},\ldots,x_m\in \X}\prod_{k=k_0+1}^m
\exp(w_{x_1,\ldots,x_k})\over
\Big(\sum_{x\in\X}\phi_v(x)\Big)\sum_{x_1,\ldots,x_m\in \X}\prod_{k=1}^m
\exp(w_{x_1,\ldots,x_k})}\\
&=c_f^m
{\Big(\prod_{k=1}^{k_0}\exp(w_{x^0_1,\ldots,x^0_k})\Big)
\sum_{x_{k+1},\ldots,x_m\in \X}\prod_{k=k_0+1}^m
\exp(w_{x_1,\ldots,x_k})\over
\sum_{x_1,\ldots,x_m\in \X}\prod_{k=1}^m
\exp(w_{x_1,\ldots,x_k})}\\
&\le c_f^m \\
&\le c_f^\Delta.
\end{align*}
Similarly,
\begin{align*}
&{\cal G}_2^{-1}{\partial {\cal G}_2\over \partial w_{x^0_1,\ldots,x^0_{k_0}}}\ge c_f^{-\Delta}.
\end{align*}
Since the dimension of the argument $\mbw$ is $\sum_{1\le k\le m}|\X|^m<\Delta \|X\|^\Delta$, then
we conclude
\begin{align*}
\|\nabla G(\mbw)\|_{L_1}\le (c_f^\Delta-c_f^{-\Delta})\Delta |X|^\Delta\le 1-\gamma.
\end{align*}
This concludes the proof.
\end{proof}

\subsection{Example: Potts model}
One of the most widely studied objects in the statistical physics is $q$-state Potts model.
It is described in the terminology of MRF as follows. Given a graph $\G$ we set $\phi_v=1$ for all nodes
$v$. $\X=\{1,2,\ldots,q\}$. A parameter $\beta$ called inverse temperature is fixed.
The coupling functions $f$ are set as $f_{u,v}(x,y)=\exp(\beta 1\{x=y\})$ for all nodes $u,v$ and all
elements $x,y\in \X$. The case $\beta>0$  corresponds to the \emph{ferromagnetic} Potts model.
In this case the distribution $\pr_{\M}(\cdot)$ "favors" assignments which select the same element
along the edges. The case $\beta<0$ corresponds to the \emph{anti-ferromagnetic} Potts model,
and in this case the distribution favors assignments with different elements along the edges.
The extreme case $\beta=-\infty$ corresponds to the usual coloring problem, where monochromatic coloring
are simply forbidden. The special case $q=2$ is called \emph{Ising} model - one of the cornerstone models
of the statistical physics.

It is immediate that conditions (\ref{eq:phinonnegative}) and (\ref{eq:fnonnegative})
are satisfied by this model provided $|\beta|<\infty$. Thus an immediate corollary of Theorem~\ref{theorem:MainResultMRF}
is the following algorithmic result.

\begin{coro}\label{coro:PottsMainResult}
There exists a deterministic FPTAS for computing a partition function for
a family of Potts model $(\G,q,\beta)$ with constant constant degree $\Delta$, constant number of colors $q$,
and satisfying
\begin{align*}
(e^{\beta\Delta}-e^{-\beta\Delta})\Delta q^\Delta< 1.
\end{align*}
\end{coro}
Observe that for large $\Delta$, the largest inverse temperature $\beta$ satisfying this condition behaves
like $O({1\over \Delta q^{\Delta}})$. We believe that this is an overly conservative estimate. We conjecture
that in fact the correlation decay property can be established in the regime
\begin{align}\label{eq:betaoverDelta}
\beta=O({1\over \Delta}),
\end{align}
leading to a deterministic FPTAS.

\section{Comparison of the correlation decay on a computation tree and the spatial correlation decay property}
\label{section:comparison}
As we have mentioned above, the (spatial) correlation decay is known to hold for the coloring problem in a stronger
regime $\alpha>\alpha^*\approx 1.763\ldots$, then the regime $\alpha>\alpha^{**}$ considered
in this paper \cite{GoldbergMartinPaterson}.
This decay of correlation is established in a conventional sense: for every node $v$ the marginal probability
$\pr(c(v)=i)$ is asymptotically independent from changing a color on a boundary of the depth-$d$ neighborhood $B(v,d)$
of $v$ in the underlying graph. In fact it is established that the decay of correlation is exponential in $d$.
It is natural to try to use this result directly as a method for computing approximately the marginals
$\pr(c(v)=i)$, for example by computing the marginal $\pr_{B(v,d)}(c(v)=i)$ corresponding to the
neighborhood $B(v,d)$, say using brute force computation.
Unfortunately, this conventional correlation decay result is not useful because of the computation
growth. In order to obtain $\epsilon$-approximation of the partition function, we need order $O(\epsilon/n)$
approximation of the marginals, which means the depth $d$ of the neighborhood $B(v,d)$ needs to be at least
$O(\log n)$. Here $n$ is the number of nodes.
But the resulting cardinality of $B(v,d)$, even for the case of constant degree graphs
is  $O(\Delta^{\log n})=n^{O(1)}$ - polynomial
in $n$ and the brute-force computation effort would be exponential in $n$. Notice that even if the underlying
graph has a polynomial expansion $|B(v,d)|\leq d^r$, for some power $r\geq 1$, the brute-force computation
would still be $O(\exp(\log^r n))$ which is super-polynomial.
This is where having correlation decay on computation tree as opposed to the conventional
graph theoretic sense helps.

\section{Conclusions}\label{section:conclusions}
We have established the existence of a deterministic approximation algorithm for
counting the number of list colorings for certain classes of graphs. We have further extended our approach
to constructing deterministic approximation algorithm for computing a partition function of a Markov random
field satisfying certain conditions.
Along with \cite{BandyopadhyayGamarnikCounting} and \cite{weitzCounting}
this work is another step in the direction of developing a new powerful method for solving counting problems
using insights from statistical physics. This method provides an important alternative to the existing
MCMC sampling based method as it leads to a deterministic as opposed to a randomized algorithm.
Since the conference version of the paper \cite{GamarnikKatz} appeared, several
new  developments happened in this direction. A deterministic approximation algorithm for counting
the number of partial matchings in constant degree graphs was constructed
in Bayati et al.~\cite{BayatiGamarnikKatzNairTetali}. The result was further used by Gamarnik
and Katz~\cite{gamarnikKatzPermanent}
for constructing a deterministic subexponential algorithm for computing a permanent of an arbitrary $0,1$ matrix.
Recently Nair and Tetali~\cite{NairTetali} introduced a somewhat different way of constructing a computation tree, closer to the
original self-avoiding walk based construction of Weitz~\cite{weitzCounting}. Furthermore, they established
that  a strong form of correlation decay (called \emph{very strong spatial mixing} in the paper) implies
correlation decay on the computation tree and ultimately leads to a polynomial time algorithm for computing
a partition function of a MRF. Their setting also allows for a hypergraph structure. It would be interesting
to use their result to perhaps tighten the condition (\ref{eq:degreealpha}) used in the present paper.

The principle insight from this work, along with the work of Weitz~\cite{weitzCounting}
is the advantage of establishing the correlation decay property on the \emph{computation
tree} as opposed to the original graph theoretic structure.
While we have established such correlation decay only in the regime $\alpha>2.8432...$, we conjecture
that it holds for much lower values of $\alpha$. In fact, just as it is conjectured that the Markov chain is rapidly
mixing in the regime $q\geq \Delta+2$, we conjecture that the correlation decay on the computation tree holds in this
regime as well, at least for the case of constant number of colors $q$. Finally, we conjecture that
the polynomial time algorithms
for computing partition function of a MRF can be constructed under weaker a assumption than (\ref{eq:conditionMRF}).
Specifically, we conjecture that for the case of Potts model, the critical inverse temperature $\beta^*$
under which the correlation decay can be established on a computation tree behaves like (\ref{eq:betaoverDelta}).

\section*{Acknowledgements} The authors are very thankful to Devavrat Shah for contributing many important comments
for this work.

\bibliographystyle{amsalpha}
\bibliography{C:/David/My_Papers/bibliography}

\end{document}